\input amstex
\documentstyle{amsppt}
\loadbold
\def\cstar{$C^*$-algebra}

\def\<{\left<}										
\def\>{\right>}

\def\tr{\text{trace}\,}

\def\esg{$E_0$-semigroup}

\nologo

\def\qedd{
\hfill
\vrule height4pt width3pt depth2pt
\vskip .5cm
}

\magnification=\magstephalf

\topmatter
\title
Interactions in noncommutative dynamics
\endtitle

\author William Arveson
\endauthor

\affil Department of Mathematics\\
University of California\\Berkeley CA 94720, USA
\endaffil

\dedicatory
To the memory of Irving Segal
\enddedicatory

\date 4 August, 1999
\enddate
\thanks On appointment as a Miller Research 
Professor in the Miller Institute for Basic 
Research in Science.  
Support is also acknowledged from 
NSF grant DMS-9802474
\endthanks
%
%
%
\abstract 
A mathematical notion of interaction is 
introduced for noncommutative dynamical systems, 
i.e., for one parameter groups of $*$-automorphisms
of $\Cal B(H)$ endowed with a certain causal structure.    
With any interaction there is a well-defined 
``state of the past" and a well-defined 
``state of the future".  
We describe the construction of many 
interactions involving cocycle perturbations
of the CAR/CCR flows and 
show that they are nontrivial.  
The proof of nontriviality is based on a new 
inequality, relating the 
eigenvalue lists of the ``past" and ``future" states
to the norm of a linear functional on 
a certain $C^*$-algebra.      
\endabstract

\rightheadtext{noncommutative interactions}

\endtopmatter

\document

\example{No minus signs?}
If you are reading this as a pdf file collected from LANL, you 
may find that it lacks 
{\bf all} minus signs.  They are missing from subscripts and 
superscripts, as well as from their normal places in formulas.  
This anomaly is mysterious to me and I don't know how to 
fix it.  
I have posted a correct pdf file for downloading from my 
Berkeley web site: 
{\tt http://www.math.berkeley.edu/\~\,arveson}.  
\endexample

\subheading{Introduction, summary of results}
In this paper we are concerned with one-parameter 
groups of $*$-automorphisms, of the algebra $\Cal B(H)$ 
of all bounded operators on a 
Hilbert space $H$, which carry a particular kind 
of causal structure.  More precisely, 
A {\bf history} is a pair $(U, M)$ consisting of 
a one-parameter group $U=\{U_t: t\in\Bbb R\}$ of unitary
operators acting on a separable infinite-dimensional
Hilbert space $H$ and a type $I$ subfactor
$M\subseteq \Cal B(H)$ which is invariant under 
the automorphisms 
$\gamma_t(X)=U_tXU_t^*$ for negative $t$, 
and which has the following two properties
\roster
\item"{(0.1)}" (irreducibility)
$$
(\bigcup_{t\in \Bbb R}\gamma_t(M))^{\prime\prime} = \Cal B(H),
$$
\item"{(0.2)}" (trivial infinitely remote past)
$$
\bigcap_{t\in\Bbb R}\gamma_t(M) = \Bbb C\cdot\bold 1.  
$$
\endroster
We find it useful to think of the group 
$\{\gamma_t: t\in\Bbb R\}$ as representing the flow
of time in the Heisenberg picture, 
and the von Neumann algebra $M$ as representing
bounded observables that are associated with the ``past".  
However, this paper is concerned with  
purely mathematical issues  
concerning the dynamical 
properties of histories, with problems concerning 
their existence and construction, and especially 
with the issue of nontriviality (to be defined momentarily).

An {\bf $E_0$-semigroup} is a one-parameter
semigroup $\alpha=\{\alpha_t: t\geq 0\}$ of unit-preserving 
$*$-endomorphisms of a type $I_\infty$ factor $M$, which 
is continuous in the natural sense 
\cite{2}--\cite{8}, \cite{10}, \cite{11}, 
\cite{29}--\cite{33}.  
The subfactors $\alpha_t(M)$ decrease as 
$t$ increases, and $\alpha$ is called 
{\bf pure} if $\cap_t\alpha_t(M)=\Bbb C\bold 1$.  
There are two $E_0$-semigroups $\alpha^-$, 
$\alpha^+$ associated with any history,
$\alpha^-$ being the one associated with 
the ``past" by restricting $\gamma_{-t}$
to $M$ for $t\geq 0$ and $\alpha^+$ being
the one associated
with the ``future" by restricting $\gamma_t$ 
to the commutant $M^\prime$ for $t\geq 0$.  
By an {\bf interaction} we mean a history with 
the additional property that there are normal
states $\omega_-$, $\omega_+$ of $M$, $M^\prime$
respectively such that $\omega_-$ is invariant 
under the action of $\alpha^-$ 
and $\omega_+$ is invariant under the action 
of $\alpha^+$.  Both $\alpha^-$
and $\alpha^+$ are pure $E_0$-semigroups, and 
when a pure $E_0$-semigroup has a normal invariant 
state then that state is uniquely determined,
see (4.1) below.  
Thus $\omega_-$ (resp. $\omega_+$) is the 
unique normal invariant state of $\alpha^-$
(resp. $\alpha^+$).

\remark{Remarks}
Since the state space of any unital \cstar\ is
weak$^*$-compact, the Markov-Kakutani fixed 
point theorem implies that every $E_0$-semigroup
has invariant states.  But there is no reason to 
expect that there is a {\it normal} invariant 
state.  Indeed, we have examples (unpublished) 
of pure $E_0$-semigroups
which have no normal invariant states.  
Notice too that $\omega_-$, for example, is defined {\it only}
on the algebra $M$ of the past.  Of course, $\omega^-$ has 
many extensions to normal states of $\Cal B(H)$, but 
none of these normal extensions 
need be invariant under the action of 
the group $\gamma$.  In fact, we will see that if there 
is a {\it normal} $\gamma$-invariant state defined 
on all of $\Cal B(H)$ then the interaction must 
be trivial.  
\endremark

In order to define a trivial interaction we must 
introduce a \cstar\ of ``local observables".  For 
every 
compact interval $[s,t]\subseteq \Bbb R$ 
there is an associated von Neumann algebra
$$
\Cal A_{[s,t]} = \gamma_t(M)\cap \gamma_s(M)^\prime.  
\tag{0.3}
$$
Notice that since $\gamma_s(M)\subseteq\gamma_t(M)$ are 
both type $I$ factors, so is the relative commutant
$\Cal A_{[s,t]}$.
Clearly $\Cal A_I\subseteq \Cal A_J$ if 
$I\subseteq J$, and for adjacent intervals 
$[r,s]$, $[s,t]$, $r\leq s\leq t$ we have 
$$
\Cal A_{[r,t]} = \Cal A_{[r,s]}\otimes\Cal A_{[s,t]},
\tag{0.4}
$$
in the sense that the two factors $\Cal A_{[r,s]}$ and 
$\Cal A_{[s,t]}$ mutually commute and generate 
$\Cal A_{[r,t]}$ as a von Neumann algebra.  The 
automorphism group $\gamma$ permutes the algebras
$\Cal A_I$ covariantly,
$$
\gamma_t(\Cal A_I) = \Cal A_{I+t}, \qquad t\in\Bbb R. 
\tag{0.5}
$$
Finally, we define the local \cstar\ $\Cal A$
to be the {\it norm} closure of the union of all 
the $\Cal A_I$, $I\subseteq \Bbb R$.  
$\Cal A$ is a $C^*$-subalgebra
of $\Cal B(H)$ which is strongly dense and 
invariant under the 
action of the automorphism group $\gamma$.  

\remark{Remarks}
It may be of interest to compare the local structure
of the \cstar\ $\Cal A$ to its commutative 
counterpart, namely the local algebras associated 
with a stationary random distribution with independent
values at every point \cite{19}.  More precisely, suppose
that we are given a random distribution $\phi$; i.e., 
a linear map from the space of real-valued test functions on 
$\Bbb R$ to the space of real-valued random variables on 
some probability space $(\Omega, P)$.  With every 
compact interval $I=[s,t]$ with $s<t$ one may consider 
the weak$^*$-closed subalgebra $\Cal A_I$ of 
$L^\infty(\Omega,P)$ generated by random variables of the form
$e^{i\phi(f)}$, $f$ ranging over all test functions supported 
in $I$.  When the random distribution $\phi$ is stationary 
and has independent values at every point, this family 
of subalgebras of $L^\infty(\Omega,P)$ 
has properties analogous to (0.4) and (0.5), in that 
there is a one-parameter group of measure preserving 
automorphisms $\gamma=\{\gamma_t: t\in \Bbb R\}$ 
of $L^\infty(\Omega,P)$ which satisfies (0.5), and 
instead of (0.4) we have the assertion that the algebras 
$\Cal A_{[r,s]}$ and $\Cal A_{[s,t]}$ are 
{\it probabilistically independent} and generate 
$\Cal A_{[r,t]}$ as a weak$^*$-closed algebra.  

One should keep in mind, however, that this commutative 
analogy has serious limitations.  For example, we have 
already pointed out that 
in the case of interactions there is 
typically no normal 
$\gamma$-invariant state on $\Cal B(H)$,
and there is no reason to expect 
any normal state of $\Cal B(H)$ to 
decompose as a product state relative to the 
decmpositions of (0.4).  

There is also some common 
ground with the Boolean algebras of type $I$ factors
of Araki and Woods \cite{1}, but here too there are 
significant differences.  For example, the local algebras
of (0.3) and (0.4) are associated with intervals (and 
more generally with finite unions of intervals), but 
not with more general Borel sets as in \cite{1}.  
Moreover, here the translation group acts as automorphisms
of the given structure whereas in \cite{1} there is 
no assumption of ``stationarity" with respect
to translations.  
\endremark

For our 
purposes, the local \cstar\ $\Cal A$ has two important features.  
First, it gives us a way of comparing $\omega_-$ 
and $\omega_+$.  
Indeed, both states $\omega_-$ and $\omega_+$ extend 
{\it uniquely} to $\gamma$-invariant states 
$\bar\omega_-$ and $\bar\omega_+$ of $\Cal A$.  
We sketch the proof for $\omega_-$.  

\proclaim{Proposition 0.6}
There is a unique $\gamma$-invariant state 
$\bar\omega_-$ of $\Cal A$ such that 
$$
\bar\omega_-\restriction_{\Cal A_I}=
\omega_-\restriction_{\Cal A_I}
$$
for every compact interval $I\subseteq (-\infty,0]$.  
\endproclaim
\demo{proof}
For existence of the extension, choose any compact interval
$I=[a,b]$ and any operator $X\in\Cal A_I$.  Then for sufficiently
large $s>0$ we have $I-s\subseteq (-\infty,0]$ and for these 
values of $s$ $\omega_-(\gamma_{-s}(X))$ does not depend on 
$s$ because $\omega_-$ is invariant under the 
action of $\{\gamma_t: t\leq 0\}$.  
Thus we can define $\bar\omega_-(X)$ unambiguously by
$$
\bar\omega_-(X)=\lim_{t\to-\infty}\omega_-(\gamma_t(X)).  
$$
This defines a positive linear functional $\bar\omega_-$ on 
the unital $*$-algebra
$\cup_I\Cal A_I$, and now we extend $\bar\omega_-$ to 
all of $\Cal A$ be norm-continuity.  The extended state is 
clearly invariant under the action of $\gamma_t$, $t\in\Bbb R$.  

The proof of uniqueness of the extension is straightforward, and we 
omit it.  
\qedd\enddemo

It is clear from the proof of Proposition 0.6 that 
these extensions of $\omega_-$ and 
$\omega_+$ are {\it locally normal} in the sense that 
their restrictions to any localized subalgebra
$\Cal A_I$ define normal states on that type $I$ factor.  

\proclaim{Definition 0.7}
The interaction $(U,M)$, with past and future states 
$\omega_-$ and $\omega_+$, is said to be trivial 
if $\bar\omega_-=\bar\omega_+$.  
\endproclaim

More generally, the norm $\|\bar\omega_- -\bar\omega_+\|$
gives some measure of the ``strength" of the 
interaction, and of course we have
$0\leq \|\bar\omega_- -\bar\omega_+\| \leq 2$.  

{\it If} there 
is a normal state $\rho$ of $\Cal B(H)$ which is invariant under 
the action of $\gamma$, then since $\omega_-$ (resp. 
$\omega_+$) is the unique normal invariant state of 
$\alpha_-$ (resp. $\alpha_+$) we must have 
$\rho\restriction_M=\omega_-$, 
$\rho\restriction_{M^\prime}=\omega_+$, and hence 
$\bar\omega_-=\bar\omega_+=\rho\restriction_\Cal A$ 
by the uniqueness part of Proposition 0.6.     
In particular, 
{\it if the interaction is nontrivial then 
neither $\bar\omega_-$ nor $\bar\omega_+$ can be extended
from $\Cal A$ to a normal state of its strong closure $\Cal B(H)$}.   

The second important feature of $\Cal A$ is that there is a definite
``state of the past" and a definite ``state of the future" 
in the following sense.  

\proclaim{Proposition 0.8}
For every $X\in \Cal A$ and every normal state $\rho$ of 
$\Cal B(H)$ we have 
$$
\lim_{t\to-\infty}\rho(\gamma_t(X))=\bar\omega_-(X), 
\qquad
\lim_{t\to+\infty}\rho(\gamma_t(X))=\bar\omega_+(X)
$$
\endproclaim
\demo{proof}
Consider the first limit formula.  
The set of all $X\in\Cal A$ for which this formula holds
is clearly closed in the operator norm, hence it suffices 
to show that it contains $\Cal A_I$ for every compact 
interval $I\subseteq\Bbb R$.  

We will make use of the fact (discussed more fully 
at the beginning of section 5) that if $\rho$ is 
any normal state of $M$ and $A$ is an operator in 
$M$ then 
$$
\lim_{t\to-\infty}\rho(\gamma_t(A))=\omega_-(A),
$$
see formula (4.1).  
Choosing a real number $T$ sufficiently negative that 
$I+T\subseteq (-\infty,0]$, the preceding remark shows
that for the operator $A=\gamma_T(X)\in M$ we have 
$\lim_{t\to-\infty}\rho(\gamma_t(A))=\omega_-(A)$, 
and hence
$$
\lim_{t\to-\infty}\rho(\gamma_t(X))=
\lim_{t\to-\infty}\rho(\gamma_{t-T}(\gamma_T(X)))=
\omega_-(\gamma_T(X))
=\bar\omega_-(X). 
$$
The proof of the second limit formula is similar.  
\qedd\enddemo

Thus, whatever (normal) state $\rho$ one chooses to watch evolve
over time on operators in $\Cal A$, it settles down to become 
$\bar\omega_+$ in the distant future, it must have 
come from $\bar\omega_-$ in the remote past, and the 
limit states do not 
depend on the choice of $\rho$.  For a trivial 
interaction, nothing happens over the long term: for 
fixed $X$ and $\rho$ 
the function $t\in\Bbb R\mapsto \rho(\gamma_t(X))$
starts out very near some value 
(namely $\bar\omega_-(X)$), exhibits transient fluctuations over 
some period of time, and then 
settles down near the same value again.  
For a nontrivial interaction, there will be a definite change 
from the limit at $-\infty$ to the limit at $+\infty$ 
(for some choices of $X\in\Cal A$).

A number of questions arise naturally.  
1) How does one construct examples of interactions?  
2) How does one determine if a given interaction 
is nontrivial? 3) What $C^*$-dynamical systems can occur as the 
$C^*$-algebras of local observables associated with an 
interaction?  
The purpose of this paper is to provide an effective partial 
solution of problem 1) and a complete solution of problem 2).   
The latter involves an inequality which we feel is of some
interest in its own right.  These results are summarized 
as follows.  

By an {\it eigenvalue list}
we mean a decreasing sequence of nonnegative
real numbers $\lambda_1\geq \lambda_2\geq\dots$ with finite 
sum.  Every normal state $\omega$ of a 
type $I$ factor is associated with a positive operator of 
trace $1$, whose eigenvalues counting multiplicity can be 
arranged into an eigenvalue list which will be 
denoted $\Lambda(\omega)$.  If the factor is finite dimensional,
we still consider $\Lambda(\omega)$ to be an infinite list by 
adjoining zeros in the obvious way.  Given two eigenvalue 
lists $\Lambda=\{\lambda_1\geq\lambda_2\geq\dots\}$ and 
$\Lambda^\prime=\{\lambda_1^\prime\geq\lambda_2^\prime\geq\dots\}$, 
we will write 
$$
\|\Lambda - \Lambda^\prime\| = 
\sum_{k=1}^\infty |\lambda_k-\lambda_k^\prime|
$$
for the $\ell^1$-distance from one list to the other.  
A classical result implies that if $\rho$ and $\sigma$ are 
normal states of a type $I$ factor $M$, then we have 
$$
\|\Lambda(\rho)-\Lambda(\sigma)\| \leq \|\rho -\sigma\| 
$$
(see section 3).  

Combining the results of \cite{7} with the results 
of section 1 below, we obtain the following 
result on the existence of interactions having 
arbitrary {\it finite} eigenvalue lists.  

\proclaim{Theorem A}
Let $n=1,2,\dots, \infty$ and let 
$\Lambda_-$ and $\Lambda_+$ be two eigenvalue lists, 
each of which has only finitely many nonzero terms.  
There is an interaction $(U,M)$ whose past and future states
$\omega_-$, $\omega_+$ have eigenvalue lists $\Lambda_-$ 
and $\Lambda_+$, and whose past and future $E_0$-semigroups
are both cocycle perturbations of the $CAR/CCR$ flow
of index $n$.   
\endproclaim

\remark{Remarks}Theorem A is established  
in section 3.  
We conjecture that the finiteness hypothesis of 
Theorem A can be dropped.  
\endremark

Theorem A gives examples of interactions, but 
it provides no information about whether or not 
these interactions are nontrivial.  
We will show that this is the case 
whenever the eigenvalue lists of $\omega_-$ and 
$\omega_+$ are different.  That conclusion depends on the 
following, which is the main result 
of this paper (and which applies 
to interactions with arbitrary...i.e., not necessarily 
finitely nonzero...eigenvalue lists).   

\proclaim{Theorem B}
Let $(U,M)$ be an interaction with past and future states
$\omega_-$ and $\omega_+$, and let $\bar\omega_-$ and 
$\bar\omega_+$ denote their extensions
to $\gamma$-invariant states of $\Cal A$.  Then 
$$
\|\bar\omega_- - \bar\omega_+\|\geq
\|\Lambda(\omega_-\otimes\omega_-) - 
\Lambda(\omega_+\otimes\omega_+)\|.  
$$
\endproclaim

\remark{Remarks}
Theorem B is proved in section 4.  
Notice the tensor product of states on the right.  For 
example, $\Lambda(\omega_-\otimes\omega_-)$ is obtained 
from the eigenvalue list 
$\Lambda(\omega_-)=\{\lambda_1\geq\lambda_2\geq\dots\}$
of $\omega_-$ by rearranging the 
doubly infinite sequence 
of all products $\lambda_i\lambda_j$, $i,j=1,2,\dots$
into decreasing order.  It can be an 
unpleasant combinatorial chore to calculate 
$\Lambda(\omega_-\otimes\omega_-)$ even when 
$\Lambda(\omega_-)$ is relatively simple and 
finitely nonzero; but we also show in 
section 4 that if $A$ and $B$ are two positive 
trace class operators such that 
$\Lambda(A\otimes A)=\Lambda(B\otimes B)$, then 
$\Lambda(A) = \Lambda(B)$.  Thus we may conclude 
\endremark

\proclaim{Corollary 1}
Let $(U,M)$, $\omega_-$, $\omega_+$ be as in Theorem B, 
and let $\Lambda_-$ and $\Lambda_+$ be the eigenvalue 
lists of $\omega_-$ and $\omega_+$ respectively.  
If $\Lambda_-\neq \Lambda_+$, then 
the interaction is nontrivial.  
\endproclaim

The following implies that ``strong" interactions exist.  

\proclaim{Corollary 2}
Let $n=1,2,\dots,\infty$ and choose $\epsilon>0$.  
There is an interaction $(U,M)$ having past and 
future states $\omega_-$, $\omega_+$, such that 
$\alpha^-$ and $\alpha^+$ are cocycle perturbations
of the $CAR/CCR$ flow 
of index $n$, for which 
$$
\|\bar\omega_- - \bar\omega_+\| \geq 2 - \epsilon.  
$$
\endproclaim

Theorem B depends on a more general result concerning
the asymptotic behavior of eigenvalue lists, 
which may be of some interest on its
own.  Let $\alpha=\{\alpha_t: t\geq 0\}$ be an 
$E_0$-semigroup acting on $\Cal B(H)$, which is 
{\it pure} in the sense defined above.  
The commutants
$N_t = \alpha_t(\Cal B(H))^\prime$ are type $I$ subfactors 
which increase with $t$, and because of 
purity their union is strongly 
dense in $\Cal B(H)$.  Let $\rho$ be a normal 
state of $\Cal B(H)$.  We require the following 
information concerning the behavior of of the 
eigenvalue lists of the restrictions 
$\rho\restriction_{N_t}$ for large $t$.  

\proclaim{Theorem C}
Let $\alpha$ be a pure $E_0$-semigroup acting on 
$\Cal B(H)$, which has a normal invariant state $\omega$.  
Then for every normal state $\rho$ of $\Cal B(H)$ we have 
$$
\lim_{t\to\infty}\|\Lambda(\rho\restriction_{\alpha_t(M)^\prime})
-\Lambda(\rho\otimes\omega)\|=0.  
$$
\endproclaim

\remark{Remarks}
One might expect that since the $N_t$ increase to 
$\Cal B(H)$, the restriction of a normal state to 
$N_t$ should look like $\rho$ itself when $t$ is large.  
Indeed, if the invariant state $\omega$ is 
a vector state then its only nonzero eigenvalue
is $1$ and $\Lambda(\rho\otimes\omega)=\Lambda(\rho)$;  
in this case Theorem C implies that 
the restriction of $\rho$ to 
$N_t$ has almost the same list as $\rho$
when $t$ is large.  On the other hand, if 
$\omega$ is not a
vector state then $\Lambda(\rho\otimes\omega)$ is 
very different from $\Lambda(\rho)$, and Theorem C
shows that this intuition is wrong.  

We also remark that Theorem C is itself a special case of a 
more general result that is independent of the theory 
of $E_0$-semigroups (see \cite{9}).  
\endremark

\subheading{1.  Existence of dynamics}

Flows on spaces are described infinitesimally by 
vector fields.  Flows on Hilbert spaces 
(that is to say, one-paramter unitary groups) are described 
infinitesimally by unbounded self-adjoint operators.  
In practice, one is usually presented with a 
symmetric operator $A$ that
is not known to be self-adjoint (much 
like being presented with a differential equation that 
is not known to posses solutions for all time), and 
one wants to know if there is a one-paramter unitary 
group that can be associated with it. Precisely, 
one wants to know if $A$ can be {\it extended} to a 
self-adjoint operator.  

This problem of the existence of dynamics was solved 
by von Neumann as follows.  Every densely defined symmetric
operator $A$ has an adjoint $A^*$ with dense domain 
$\Cal D^*$, and using $A^*$ one defines two 
{\it deficiency spaces} $\Cal E_-$, $\Cal E_+$ by
$$
\Cal E_{\pm} =\{\xi\in\Cal D^*: A^*\xi=\pm i\xi\}.  
$$ 
von Neumann's result is that $A$ has self-adjoint 
extensions iff $\dim\Cal E_-=\dim\Cal E_+$ (see 
\cite{15, section XII.4}).  Moreover, when 
$\Cal E_-$ and $\Cal E_+$ have the same dimension, 
von Neumann showed that for every unitary operator
from $\Cal E_-$ to $\Cal E_+$ there is an 
associated self-adjoint extension of $A$.  
The purpose of this section is to establish an 
analogous result which locates the obstruction 
to the existence of dynamics for pairs of 
$E_0$-semigroups of the simplest kind 
(Corollary 1 below).  That is based on the 
following more general result.  

Let $M$ be a type $I$ subfactor of $\Cal B(H)$, and 
let $\alpha$, $\beta$ be two $E_0$-semigroups acting, 
respectively, on $M$ and its commutant $M^\prime$.  
We want to examine conditions under which there is a 
one-parameter unitary group $U=\{U_t: t\in\Bbb R\}$ 
acting on $H$ whose associated automorphism group
$\gamma_t(A)=U_tAU_t^*$ has $\alpha$ as its past
and $\beta$ as its future in the sense 
that 
$$
\gamma_{-t}\restriction_{M}=\alpha_t,\qquad
\gamma_t\restriction_{M^\prime}=\beta_t,
\qquad t\geq 0.
\tag{1.1}
$$
The following result asserts that there is such a 
unitary group $U$ if and only if the 
product systems of $\alpha$ and $\beta$ 
are {\it anti-}isomorphic.  

\proclaim{Theorem}
Let $E^\alpha=\{E^\alpha(t): t>0\}$ 
and $E^\beta=\{E^\beta(t): t>0\}$ be the 
respective product systems of $\alpha$ and $\beta$, 
$$
\align
E^\alpha(t)&=\{x\in M: \alpha_t(y)x=xy, \qquad y\in M\},\\
E^\beta(t)&=\{x^\prime\in M^\prime: \beta_t(y^\prime)x^\prime
=x^\prime y^\prime, \quad y^\prime\in M^\prime\},  
\endalign
$$
and assume that 
there is a one-parameter 
unitary group $U=\{U_t: t\in\Bbb R\}$ whose associated 
automorphism group satisfies (1.1).  Then 
$E^\alpha$ and $E^\beta$ are anti-isomorphic.  Indeed, 
for every $t>0$ we have $U_tE^\alpha(t) = E^\beta(t)$, 
and the map $\theta: E^\alpha\to E^\beta$ defined 
by 
$$
\theta(v) = U_tv,\qquad v\in E^\alpha(t), \quad t>0,
\tag{1.2}
$$
is an anti-isomorphism of product systems (i.e., it is 
a Borel-measurable map which is unitary on fibers, and 
which satisfies $\theta(vw)=\theta(w)\theta(v)$ for 
every $v\in E^\alpha(s)$, $w\in E^\alpha(t)$, 
$s,t >0$).  

Conversely, if $\theta: E^\alpha\to E^\beta$ is any 
anti-isomorphism of product systems, then for every 
$t>0$ there is a unique unitary operator 
$U_t\in\Cal B(H)$ which satisfies (1.2) for every 
$v\in E^\alpha(t)$.  $\{U_t: t>0\}$ is a strongly 
continuous semigroup of unitary operators tending
strongly to the identity as $t\to 0+$, 
and its natural extension to a one-parameter unitary 
group gives rise to an automorphism group $\gamma$ which 
satisfies (1.1).  
\endproclaim

\demo{proof}
Assume that $\gamma_t(A)=U_tAU_t^*$, $t\in \Bbb R$ 
satisfies (1.1).  Fix $t>0$.  We claim first that 
$U_tE^\alpha(t)\subseteq M^\prime$.  
Indeed, if $x\in M$ then for every $v\in E^\alpha(t)$
we have 
$$
xU_tv = U_t\gamma_{-t}(x)v = U_t\alpha_t(x)v = U_tvx.  
$$

Next, we claim that $U_tE^\alpha(t)\subseteq E^\beta(t)$.  
For $v\in E^\alpha(t)$, the preceding shows that 
$U_tv\in M^\prime$, so it suffices to show that 
$\beta_t(y)U_tv = U_tvy$ for every $y\in M^\prime$.  For 
that, write
$$
\beta_t(y)U_tv = \gamma_t(y)U_tv =
U_tyU_t^*U_tv = U_tyv = U_tvy,
$$
the last equality because $v\in M$ commutes with $y\in M^\prime$.

Next, note that $E^\beta(t)\subseteq U_tE^\alpha(t)$.  
Choosing $w\in E^\beta(t)$, set $v=U_t^*w$.  Note that 
$v\in M$ because for every $y\in M^\prime$ we have 
$$
yv = yU_t^*w= U_t^*\gamma_t(y)w=U_t^*\beta_t(y)w =
U_t^*wy = vy.  
$$

Note next that the element $v=U_t^*w\in M$ actually belongs
to $E^\alpha(t)$.  Indeed, for every $x\in M$ we have 
$$
\alpha_t(x)v=\alpha_t(x)U_t^*w.  
$$
Since $\gamma_{-t}$ restricts to $\alpha_t$ on $M$, 
we have $\gamma_t(\alpha_t(x)) = x$ and 
the right side can be written 
$$
U_t^*\gamma_t(\alpha_t(x)) = U_t^*xw = U_t^*wx = vx.  
$$

The above shows that for every $t>0$ we have a linear 
map $\theta_t: E^\alpha(t)\to E^\beta(t)$ defined by 
$\theta_t(v)=U_tv$.  By assembling these maps we get a 
Borel-measurable map $\theta: E^\alpha\to E^\beta$ 
which is linear on fibers.  Notice that $\theta_t$ 
is actually unitary, since for $v_1,v_2\in E^\alpha(t)$
we have 
$$
\<v_1,v_2\>\bold 1 = v_2^*v_1= 
(U_tv_2)^*(U_tv_1)=\theta(v_2)^*\theta(v_1)=
\<\theta(v_1),\theta(v_2)\>\bold 1.  
$$
Finally, $\theta$ is an anti-isomorphism, because 
for $v\in E^\alpha(s)$, $w\in E^\alpha(t)$ we have 
$$
\theta(vw) = U_{s+t}vw=U_t(U_sv)w=U_t\theta(v)w =
U_tw\theta(v)=\theta(w)\theta(v).  
$$

To prove the converse, fix an anti-isomorphism 
$\theta: E^\alpha\to E^\beta$.  For every $t>0$ 
pick an orthonormal basis $e_1(t),e_2(t)\dots$ for 
$E^\alpha(t)$ (we will have to choose more carefully 
presently...but for the 
moment we choose an arbitrary orthonormal
basis for each fiber space).  For every $t>0$ define 
an operator $U_t\in \Cal B(H)$ by 
$$
U_t = \sum_{n=1}^\infty \theta(e_n(t))e_n(t)^*.  
$$
One checks easily that $U_tU_t^*=U_t^*U_t=\bold 1$, hence
$U_t$ is unitary.  $U_t$ also satisfies (1.2), 
for if $v\in E^\alpha(t)$ then 
we have $e_n(t)^*v=\<v,e_n(t)\>\bold 1$ and hence 
$$
U_tv = \sum_{n=1}^\infty \<v,e_n(t)\>\theta(e_n(t)) =
\theta(\sum_n\<v,e_n(t)\>e_n(t))=\theta(v).  
$$
Note too that since the ranges of the operators in 
$E^\alpha(t)$ span $H$, any operator $U_t$ that satisfies
(1.2) is determined uniquely.  In particular, $U_t$ does not
depend on the choice of orthonormal basis $\{e_n(t)\}$ 
for $E^\alpha(t)$.  

We may choose the orthonormal basis $\{e_n(t)\}$ so that 
each section $t\mapsto e_n(t)\in E^\alpha(t)$ is Borel 
measurable (because of the measurability axiom of product 
systems \cite{2, Property 1.8 (iii)}), 
and once this is done we find that the function 
$t\in (0,\infty)\mapsto U_t\in \Cal B(H)$ is Borel measurable.  

We claim next that $\{U_t: t>0\}$ is a semigroup.  Indeed, 
if $w\in E^\alpha(s)$, $v\in E^\alpha(t)$ then since 
$\theta(v)\in M^\prime$ commutes with $w\in M$ we have 
$$
U_sU_tvw=U_s\theta(v)w=U_sw\theta(v)=\theta(w)\theta(v) =
\theta(vw)=U_{s+t}vw.  
$$
Since $E^\alpha(s+t)$ is spanned by such product $vw$ and 
since $E^\alpha(s+t)H$ spans $H$, we conclude that 
$U_sU_t=U_{s+t}$.  

At this point, we use the measurability proposition 
\cite{2, Proposition 2.5 (ii)} 
(stated there for the more general case of 
cocycles) to conclude that a) $U_t$ is strongly continuous
in $t$ for $t>0$, and b) $U_t$ tends strongly to $\bold 1$ 
as $t\to 0+$.  Now extend $U$ in the obvious way to obtain a 
strongly continuous one-parameter unitary group acting on 
$H$.  

Let $\gamma_t(A) = U_tAU_t^*$, $A\in\Cal B(H)$, $t\in\Bbb R$.  
It remains to show that for every $t>0$ we have 
$\gamma_{-t}\restriction_M=\alpha_t$ and 
$\gamma_t\restriction_{M^\prime}=\beta_t$.  

Choose $x\in M$.  To show that $\gamma_{-t}(x)=\alpha_t(x)$,
it suffices to show that $\gamma_{-t}(x)v=\alpha_t(x)v$
for every $v\in E^\alpha(t)$ (because $H$ is spanned by the 
ranges of the operators in $E^\alpha(t)$).  But for such a
$v$ we have 
$$
\gamma_{-t}(x)v=U_{-t}xU_tv=U_{-t}x\theta(v)=
U_{-t}\theta(v)x=vx=\alpha_t(x)v.  
$$

Choose $y\in M^\prime$.  To show that $\gamma_t(y)=\beta_t(y)$
it suffices to show that $\gamma_t(y)w=\beta_t(y)w$ 
for all $w\in E^\beta(t)$.  For such a $w$ we have 
$w=\theta(v)=U_tv$ for some $v\in E^\alpha(t)$, hence 
$$
\gamma_t(y)w=U_tyU_t^*U_tv = U_tyv=U_tvy=wy=\beta_t(y)w,
$$
and the proof is complete\qedd
\enddemo

We view the following result as a counterpart for noncommutative 
dynamics of von Neumann's theorem on the existence of self-adjoint
extensions of symmetric operators in terms of deficiency indices.   

\proclaim{Corollary 1}
Let $\alpha$ and $\beta$ be two $E_0$-semigroups, acting 
on $\Cal B(H)$ and $\Cal B(K)$ respectively, 
each of which is a cocycle perturbation of 
a CCR/CAR flow.  There is a one-parameter group of automorphisms
of $\Cal B(H\otimes K)$ which satisfies the condition
of (1.1) if, and only if, $\alpha$ and $\beta$ have the same
numerical index.  
\endproclaim

\demo{proof}
Consider the type $I$ subfactor $M$ 
of $\Cal B(H\otimes K)$ defined by
$$
M=\Cal B(H)\otimes\bold 1_K.  
$$
We have $M^\prime=\bold 1_H\otimes\Cal B(K)$, and 
$\alpha$ (resp. $\beta$) is conjugate to the 
action on $M$ (resp. $M^\prime$) defined 
by $A\otimes\bold 1_K\mapsto \alpha_t(A)\otimes\bold 1_K$
(resp. $\bold 1_H\otimes B\mapsto \bold1_H\otimes\beta_t(B)$),
$t\geq 0$.  

Now the product system of 
any $CAR/CCR$ flow is anti-isomorphic
to itself.  This follows, for example, from the structural 
results on divisible product systems of \cite{2, section 6}.
Alternately, one can simply write down 
explicit anti-automorphisms of the product systems 
described on pp. 12--14 of \cite{2}.  Since the structure
of the product system of any $E_0$-semigroup is stable 
under cocycle perturbations, the same is true of 
cocycle perturbations of CAR/CCR flows.  

The preceding theorem implies that there is a one-parameter
group of automorphisms $\gamma=\{\gamma_t: t\in\Bbb R\}$ of 
$\Cal B(H\otimes K)$ satisfying 
$$
\gamma_{-t}(A\otimes\bold 1_K)=\alpha_t(A)\otimes \bold 1_K, \qquad
\gamma_{t}(\bold 1_H\otimes B)=\bold 1_H\otimes \beta_t(B)
$$
for every $t\geq 0$ iff the product systems $E^\alpha$ and 
$E_\beta$ are anti-isomorphic. The preceding paragraph
shows that this is true iff $E_\alpha$ and $E_\beta$ are 
isomorphic; and since $\alpha$ and $\beta$ are simply cocycle 
perturbations of CAR/CCR flows, the latter holds iff 
$\alpha$ and $\beta$ have the same numerical index.  \qedd\enddemo

\proclaim{Corollary 2}
Let $\alpha$ and $\beta$ be two pure $E_0$-semigroups
which are cocycle-conjugate to the $CAR/CCR$ flow
of index $n=1,2,\dots,\infty$.  Then there is a history
$(U,M)$ whose past and future semigroups are conjugate, 
respectively, to $\alpha$ and $\beta$.  
\endproclaim

\subheading{2.  Eigenvalue lists of normal states}

In this section we emphasize the importance
of the ``eigenvalue list" invariant that 
can be associated with 
normal states of type $I$ factors, and we 
summarize its basic properties.   
An {\it eigenvalue list} is a decreasing sequence 
$\lambda_1\geq\lambda_2\geq\dots$ of nonnegative 
real numbers satisfying $\sum_n\lambda_n<\infty$.  If 
$\Lambda=\{\lambda_1\geq\lambda_2\geq\dots\}$ and 
$\Lambda^\prime=
\{\lambda_1^\prime\geq\lambda_2^\prime\geq\dots\}$
are two such lists we write 
$$
\|\Lambda - \Lambda^\prime\|=\sum_{n=1}^\infty 
|\lambda_n-\lambda_n^\prime|
$$
for the $\ell^1$-distance from $\Lambda$ to
$\Lambda^\prime$, thereby making the space of 
all eigenvalue lists into a complete metric space.  

Let $A$ be a positive trace class operator acting on 
a separable Hilbert space $H$.  
The positive eigenvalues of $A$ (counting multiplicity)
can be arranged in decreasing order, and if there are only 
finitely many nonzero eigenvalues then we extend the list 
by appending zeros in the obvious way.  This defines the 
eigenvalue list $\Lambda(A)$ of $A$.  Notice that even 
when $H$ is finite dimensional, $\Lambda(A)$ is an infinite
list.  

The following basic properties of eigenvalue lists 
will be used repeatedly.  

\proclaim{Proposition 2.1}
\roster
\item"{2.1.1}"
For every positive trace class operator $A$ we have 
$\Lambda(A)=\Lambda(A\oplus 0_\infty)$, $0_\infty$ 
denoting the infinite dimensional zero operator.  
\item"{2.1.2}"
For positive trace class operators $A$ and $B$, 
$\Lambda(A)=\Lambda(B)$ iff $A\oplus 0_\infty$ is 
unitarily equivalent to $B\oplus0_\infty$.  
\item"{2.1.3}"
If $L$ is any Hilbert-Schmidt operator from a Hilbert space 
$H_1$ to a Hilbert space $H_2$, then 
$\Lambda(L^*L)=\Lambda(LL^*)$.  
\item"{2.1.4}"
For positive trace class operators $A$, $B$ we have 
$\Lambda(A)=\Lambda(B)$ iff 
$$
{\text{trace}}(A^n)={\text{trace}}(B^n)\qquad 
{\text{for every }}n=1,2,\dots.  
$$
\endroster
\endproclaim
\demo{proof}
The assertion (2.1.1) is obvious, and (2.1.2) follows after a 
routine application of the spectral theorem for self-adjoint 
compact operators.  

\demo{proof of (2.1.3)}
Let $K_1\subseteq H_1$ be the initial space of $L$ and 
let $K_2=\overline{LK_1}\subseteq H_2$ be its closed range.   
The polar decomposition implies that $L^*L\restriction_{K_1}$ 
and $LL^*\restriction_{K_2}$ are unitarily equivalent.  
Hence $L^*L\oplus0_\infty$ and $LL^*\oplus_\infty$ are 
unitarily equivalent and the assertion 
(2.1.3) follows from (2.1.2).  
\enddemo

\demo{proof of (2.1.4)}
If $\Lambda(A)=\{\lambda_1\geq\lambda_2\geq\dots\}$ then 
$$
{\text{trace}}(A^n)=\sum_{k=1}^\infty \lambda_k^n, 
\qquad n=1,2,\dots.  
$$
Thus $\Lambda(A)=\Lambda(B)$ implies that 
${\text{trace}}(A^n)={\text{trace}}(B^n)$ for every $n\geq 1$.  

Conversely, suppose that ${\text{trace}}(A^n)={\text{trace}}(B^n)$
for every $n=1,2,\dots$.  Choose a positive number $M$ so large that 
the interval $[0,M]$ contains the spectra of both 
operators $A$ and $B$.  
The linear functional $f\mapsto {\text{trace}}(Af(A))$ defined
on the commutative \cstar\ $C[0,M]$ is positive, hence there is a 
unique finite positive measure 
$\mu_A$ defined on $[0,M]$ such that 
$$
\int_0^M f(x)\,d\mu_A(x)={\text{trace}}(Af(A)), \qquad
f\in C[0,M].  
$$
The restriction of $\mu_A$ to $(0,M]$ is concentrated on 
$\sigma(A)\cap (0,M]$, and for every positive eigenvalue $\lambda$ 
of $A$ we have 
$$
\mu_A(\{\lambda\})=\lambda\cdot {\text{multiplicity of }}\lambda.
$$

Doing the same for the operator $B$, we find that by hypothesis 
$$
\int_0^M x^n\,d\mu_A(x)=\int_0^M x^n\,d\mu_B(x), \qquad n=0,1,2,\dots,
$$
and hence by the Weierstrass approximation theorem $\mu_A$ 
and $\mu_B$ define the same linear functional on $C[0,M]$.  
It follows that $\mu_A=\mu_B$,  and the preceding observations 
lead us to conclude that $\Lambda(A)=\Lambda(B)$.  
\qedd
\enddemo
\enddemo

We will also make use of the following classical result, 
originating in work of Hermann Weyl around 1912.  

\proclaim{Proposition 2.2}
If $A$, $B$ are positive trace class operators acting on 
the same Hilbert space $H$, then 
$$
\|\Lambda(A)-\Lambda(B)\|\leq {\text{trace}}|A-B|.  
$$
\endproclaim
\demo{proof}
A proof can be found in the appendix of \cite{29}.  
\qedd\enddemo

\remark{Remarks}
Notice that since $\Lambda(A)$ depends only on the unitary 
equivalence class of $A$, Proposition 2.2 actually implies that 
$$
\|\Lambda(A)-\Lambda(B)\| \leq \inf_{A^\prime, B\prime}
{\text{trace}}|A^\prime-B^\prime|,
$$
where $A^\prime$ (resp. $B^\prime$) ranges over all operators
unitarily equivalent to $A$ (resp. $B$).  Indeed, 
though we do not require the fact, it is not hard to 
show that $\|\Lambda(A)-\Lambda(B)\|$ is exactly the 
distance (relative to the trace norm) from the unitary 
equivalence class of $A\oplus0_\infty$ to the unitary 
equivalence class of $B\oplus0_\infty$.  Thus the eigenvalue
list $\Lambda(A)$ provides a 
more-or-less complete invariant for classifying positive 
trace class operators up to unitary equivalence.  

On the other hand, 
the eigenvalue list is also a subtle invariant.  
To illustrate the point, suppose that $A$ has only two positive 
eigenvalues $3/4$ and $1/4$, and that 
$B$ has only three positive eigenvalues $3/5,1/5,1/5$.    
The spectrum of $A\oplus B$ is the union of the spectra 
and the spectrum of $A\otimes B$ is the set of products of elements
from the two spectra; however, both of these sets must 
be rearranged in decreasing order.  Thus
$$
\align
\Lambda(A\oplus B)&=\{3/4, 3/5, 1/4,1/5,1/5,0,\dots\},\\
\Lambda(A\otimes B)&=\{9/20,3/20,3/20,3/20,1/20,1/20,0,\dots\}.  
\endalign
$$
Notice that $A$ has only eigenvalues of multiplicity $1$, 
$B$ has eigenvalues of multiplicities $1$ and $2$, 
but that $A\otimes B$ has an eigenvalue of ``peculiar" 
multiplicity $3$.
In the case of larger spectra, the relation between 
say $\Lambda(A\otimes B)$ and the individual
lists $\Lambda(A)$ and $\Lambda(B)$ depends in a 
complex way on the relative sizes of eigenvalues, and 
the problem of rearranging the set of products into 
decreasing order can be a difficult combinatorial chore.  
\endremark

Turning now to normal states, let $M$ be a type $I_n$ factor,
$n=1,2,\dots,\infty$ (one can assume without essential loss that $M$
is concretely representated as a subfactor of $\Cal B(H)$ for 
some Hilbert space $H$), and let $\rho$ be a normal state 
of $M$.  There is a Hilbert space $K$ of dimension $n$ such 
that $M$ is isomorphic as a $*$-algebra to $\Cal B(K)$, and 
in this case any such $*$-isomorphism must be isometric and 
normal.  Thus we may identify $\rho$ with a normal state 
of $\Cal B(K)$, and consequently there is a positive operator 
$R\in\Cal B(K)$ of trace $1$ such that 
$$
\rho(T)={\text{trace}}(RT),\qquad T\in\Cal B(K).  
$$
The {\it eigenvalue list} of $\rho$ is defined by 
$\Lambda(\rho)=\Lambda(R).$
The preceding discussion leads immediately to the following. 

\proclaim{Proposition 2.3}
\roster
\item"{2.3.1}"
If $\rho_1$ and $\rho_2$ are normal states of type $I$ factors
$M_1$ and $M_2$, and if $\rho_1$ and $\rho_2$ are 
conjugate in the sense that there is a $*$-isomorphism 
$\theta$ of $M_1$ onto $M_2$ 
such that $\rho_2\circ\theta=\rho_1$, 
then $\Lambda(\rho_1) = \Lambda(\rho_2)$.  
\item"{2.3.2}"
If $\rho_1$ and $\rho_2$ are two normal states of a type 
$I$ factor $M$, then 
$$
\|\Lambda(\rho_1)-\Lambda(\rho_2)\|\leq \|\rho_1-\rho_2\|.
$$
\endroster
\endproclaim
\demo{proof}
The first assertion is apparent after we we realize 
$M_k$ as $\Cal B(H_k)$, $k=1,2$, use the fact that 
a $*$-isomorphism of $\Cal B(H_1)$ onto $\Cal B(H_2)$ is 
implemented by a unitary operator from $H_1$ to $H_2$, and 
make use of (2.1.2).  
The second assertion is the inequality of Proposition 2.2.  
\qedd\enddemo

\subheading{3.  CP semigroups and the existence of interactions}

The corollary of section 1 implies that any pair of 
pure $E_0$-semigroups $\alpha_-$, $\alpha_+$, which 
are both cocycle conjugate to the same $CAR/CCR$
flow, can be assembled so as to obtain a history 
$(U,M)$ whose past and future $E_0$-semigroups are 
conjugate to 
$\alpha_-$ and $\alpha_+$.  Moreover, if both 
$\alpha_-$ and $\alpha_+$ have normal invariant 
states then $(U,M)$ is in fact an interaction.  

Thus we are led to ask what the possibilities 
are.  More precisely, suppose we are given an 
eigenvalue list 
$\Lambda=\{\lambda_1\geq\lambda_2\geq\dots\}$
with $\sum_n\lambda_n=1$ and a nonnegative 
integer $n=1,2,\dots,\infty$.  Does there exist a 
cocycle perturbation $\alpha$ of the $CAR/CCR$ flow
of index $n$ which is pure, and which leaves invariant
a normal state whose eigenvalue list is $\Lambda$?  

We do not know the answer in general, but we 
conjecture that it is yes.  The purpose of this section 
is to provide an affirmative answer for the cases in 
which $\Lambda$ has only a finite number of nonzero terms
(Theorem A).  
This is essentially the main result of \cite{7} (together 
with Corollary 1 of section 1), and 
we merely summarize the main ideas so 
as to emphasize the 
role of dilation theory and semigroups of 
completely postive maps (sometimes called quantum dynamical 
semigroups) acting on matrix algebras, for 
such constructions. 

Suppose that $\alpha=\{\alpha_t: t\geq 0\}$ is 
an $E_0$-semigroup acting on $\Cal B(H)$, and 
assume further that there is a normal state 
$\omega$ of $\Cal B(H)$ which is invariant,
$\omega\circ\alpha_t=\omega$, $t\geq 0$.  Letting 
$\Omega$ be the density operator of $\omega$,
$$
\omega(T)={\text{trace}}(\Omega T), \qquad T\in\Cal B(H)
$$
then the projection $P$ on the closed range of $\Omega$ 
is the support projection of $\omega$, i.e., the largest 
projection with the property that $\omega(P^\perp)=0$.  
Using $\omega\circ\alpha_t=\omega$, we find that 
$\omega(\bold 1-\alpha_t(P))=\omega(\alpha_t(P^\perp))=
\omega(P^\perp)=0$, hence 
$\bold 1-\alpha_t(P)\leq\bold 1- P$, hence 
$$
\alpha_t(P)\geq P,\qquad t\geq 0.  
\tag{3.1}
$$ 

The inequality (3.1) has the 
following consequence.  If we identify
$\Cal B(PH)$ with the corner 
$P\Cal B(H)P$, then for every 
$t\geq 0$ we can compress $\alpha_t$ so as to obtain a 
completely positive map $\phi_t$ on $\Cal B(PH)$ 
$$
\phi_t(X) = P\alpha_t(X)\restriction_{PH}, 
\qquad X\in P\Cal B(H)P.  
$$ 
More significantly, because of (3.1) we have the semigroup
property $\phi_s\circ\phi_t=\phi_{s+t}$, as one can easily 
verify using $P\alpha_s(A)P=P\alpha_s(PAP)P$ for $A\in\Cal B(H)$.  
Thus we have defined a semigroup $\phi=\{\phi_t: t\geq 0\}$ 
of normal completely positive maps of $\Cal B(PH)$ satisfying 
$\phi_t(\bold 1)=\bold 1$ for $t\geq 0$, together with the 
natural continuity property 
$$
\lim_{t\to t_0}\<\phi_t(X)\xi,\eta\>=\<\phi_{t_0}(X)\xi,\eta\>,
$$
$\xi,\eta\in PH$, $X\in\Cal B(PH)$.  

We appear to have lost ground, in that we started with a 
semigroup of $*$-endomorphisms and now have merely a 
semigroup of completely positive maps.  However, notice 
that the restriction of $\omega$ to $\Cal B(PH)=P\Cal B(H)P$ 
is a {\it faithful} normal state which is invariant under the 
action of $\phi$, $\omega\circ\phi_t=\omega$, $t\geq 0$.  

Notice too that in case there are only a finite number of 
positive eigenvalues in the list $\Lambda(\omega)$ then 
$PH$ is finite dimensional, and thus $\phi=\{\phi_t: t\geq 0\}$
is a CP semigroup acting essentially on a {\it matrix algebra}, 
which leaves invariant a faithful state with prescribed 
eigenvalues $\lambda_1\geq\lambda_2\geq\dots\geq\lambda_r>0$.  
If $\alpha$ began life as a pure $E_0$-semigroup then 
$\omega$ is an {\it absorbing} state for $\phi$ in 
the sense that for every normal state $\rho$ of $\Cal B(PH)$
$$
\lim_{t\to\infty}\|\rho\circ\phi_t-\omega\|=0.  
\tag{3.2}
$$

Conversely and most significantly, if we can 
create a pair $(\phi,\omega)$ satisfying 
the conditions of the preceding paragraph then it
is possible to reconstruct
a pair $(\alpha,\omega)$ consisting of an $E_0$-semigroup
$\alpha$ having an invariant normal state $\omega$ with the expected 
eigenvalue list by a ``dilation" procedure which 
reverses the ``compression" procedure we have described
above.  Moreover, 
if the CP semigroup $\phi$ has a bounded generator (as it 
will surely have in the case where $PH$ is finite dimensional), 
then its dilation to an $E_0$-semigroup will be cocycle-conjugate
to a $CAR/CCR$ flow whose index can be calculated directly in 
terms of $\phi$ (the details can be found in
\cite{7} and \cite{8}).  The following summarizes the result of 
the construction of $(\phi,\omega)$ for finite eigenvalue lists 
given in \cite{7, Theorem 5.1}.

\proclaim{Theorem 3.3}
Let $\lambda_1\geq\lambda_2\geq\dots\geq\lambda_r>0$ 
be a list of positive numbers and let $\omega$ 
be a state of the matrix algebra $M_r(\Bbb C)$ whose 
density operator has this ordered list of eigenvalues.  

There is a semigroup $\phi=\{\phi_t: t\geq 0\}$ 
of unital completely positive maps on $M_r(\Bbb C)$ 
which leaves $\omega$ invariant, satisfies (3.2), and 
which can be dilated to a pure cocycle perturbation of 
a $CAR/CCR$ flow having a normal invariant state 
whose eigenvalue list has exactly 
$\lambda_1\geq\dots\geq\lambda_r$ as its nonzero elements.
\endproclaim

Theorem 3.3 leads to the following (see pp. 40--42 of \cite{7}).  

\proclaim{Corollary}
Let $n=1,2,\dots,\infty$ and let 
$\Lambda=\{\lambda_1\geq\lambda_2\geq\dots\}$ be an 
eigenvalue list which has only a finite number of 
nonzero terms.  There is a cocycle perturbation 
$\alpha$ of the $CAR/CCR$ flow of index $n$ which is 
pure, and which has an invariant normal state with 
eigenvalue list $\Lambda$.  
\endproclaim

Using Corollary 1 of section 1, we  
deduce Theorem A of the introduction.  

\proclaim{Theorem A}
Let $n=1,2,\dots,\infty$ and let $\Lambda_-$ and 
$\Lambda_+$ be two eigenvalue lists having only 
a finite number of nonzero terms.  There is an 
interaction $(U,M)$ whose past and future normal 
states $\omega_-$, $\omega_+$ have eigenvalue 
lists $\Lambda_-$, $\Lambda_+$ respectively, 
and whose past and future $E_0$-semigroups 
are cocycle conjugate to the $CAR/CCR$ flow 
of index $n$.  
\endproclaim

\subheading{4.  The interaction inequality}

Theorem A provides many examples of interactions, 
but it says nothing about whether or not these 
interactions are nontrivial.  For that we need 
the inequality of Theorem B of the introduction.  
The purpose of this section is to prove  
Theorem B and discuss 
its consequences for interactions.  Theorem B is 
based on the following more general 
result about \esg s.  An 
\esg\ $\alpha=\{\alpha_t: t\geq 0\}$ acting on 
$\Cal B(H)$ is said to be {\it pure} if 
$$
\bigcap_{t\geq 0}\alpha_t(\Cal B(H)) = \Bbb C\cdot\bold 1.  
$$
Purity implies that for any two normal states 
$\rho_1$ and $\rho_2$
$$
\lim_{t\to\infty}\|\rho_1\circ\alpha_t-\rho_2\circ\alpha_t\|=0
$$
see Proposition 1.1 of \cite{7}.  In particular, if there is 
a {\it normal} state $\omega$ which is invariant 
under $\alpha$ in the sense that 
$\omega\circ\alpha_t=\omega$ for every $t\geq 0$ 
then $\omega$ must be an {\it absorbing} state in the 
sense that for every normal state $\rho$ of 
$\Cal B(H)$ we have 
$$
\lim_{t\to\infty}\|\rho\circ\alpha_t-\omega\|=0.  
\tag{4.1}
$$
Thus, if a pure \esg\ has a normal invariant
state then it is unique, and in particular the eigenvalue list 
$\Lambda(\omega)$ of a normal invariant state $\omega$ provides 
a conjugacy invariant of pure \esg s.

Given a pure \esg\ acting on $\Cal B(H)$, the commutants
$N_t=\alpha_t(\Cal B(H))^\prime$ are type $I$ subfactors of 
$\Cal B(H)$ which increase with $t$, and 
by purity their union
is a strongly dense $*$-subalgebra of $\Cal B(H)$.  
Let $\rho$ be any normal state of $\Cal B(H)$.  Since 
$N_t$ is a type $I$ factor, the restriction of 
$\rho$ to $N_t$ has an eigenvalue list, defined as 
in section 3.  
The following result shows 
how these eigenvalue lists behave for large $t$.  

\proclaim{Theorem C}
Let $\alpha=\{\alpha_t: t\geq 0\}$ be a pure 
\esg\ having a normal invariant state $\omega$, 
and let $N_t$ be the commutant 
$\alpha_t(\Cal B(H))^\prime$.  
Then for every normal state $\rho$ of $\Cal B(H)$ we have 
$$
\lim_{t\to\infty}\|\Lambda(\rho\restriction_{N_t})-
\Lambda(\rho\otimes\omega)\|=0.  
$$
\endproclaim

The proof of Theorem C requires some preparation.  

\proclaim{Lemma 4.2}
Let $\{A_i: i\in I\}$ be a net of positive trace 
class operators acting on a Hilbert space $H$ and 
let $B$ be a positive trace class operator such 
that $\tr (A_i)=\tr (B)$ for every $i\in I$.  Suppose
there is a set $S\subseteq H$, having $H$ as its closed 
linear span, such that 
$$
\lim_i \<A_i\xi,\eta\>=\<B\xi,\eta\>, \qquad \xi,\eta\in S.  
$$
Then ${\text{trace}} |A_i-B|\to 0$, as $i\to\infty$.  
\endproclaim
\demo{proof}
By Proposition 1.6 of \cite{7} it suffices to 
show that 
$$
\lim_{i\to\infty}{\text{trace}}(A_iK)={\text{trace}}(BK)
$$
for every compact operator $K\in\Cal B(H)$.  The set 
$\Cal S$ of  compact operators $K$ for which the assertion 
is true is a norm-closed linear space which contains 
all rank-one operators of the form
$\zeta\mapsto \<\zeta,\xi\>\eta$, with $\xi,\eta\in S$.  
Since $S$ spans $H$, it follows
that $\Cal S$ is the space of all compact operators.  
\qedd
\enddemo

The next three Lemmas relate to the following situation.  
We are given a normal $*$-endomorphism $\alpha$ of 
$\Cal B(H)$ satisfying $\alpha(\bold 1)=\bold 1$.  Let 
$\Cal E$ be the linear space of operators
$$
\Cal E=\{v\in\Cal B(H): \alpha(x)v=vx, \quad x\in\Cal B(H)\}.  
$$
If $u,v$ are any two elements of $\Cal E$ then $v^*u$ is 
a scalar multiple of the identity operator, and in fact 
$\Cal E$ is a Hilbert space relative to the inner product 
defined on it by 
$$
v^*u = \<u,v\>_{\Cal E}\bold 1.  
$$
For any orthonormal basis $v_1, v_2,\dots$ of $\Cal E$ 
we have 
$$
\alpha(x) = \sum_n v_nxv_n^*, \qquad x\in\Cal B(H).  
$$

Let $\rho$ be a normal state of $\Cal B(H)$.  It is clear 
that $u,v\in\Cal E\mapsto \rho(uv^*)$ defines a bounded 
sesquilinear form on the Hilbert space $\Cal E$, hence by 
the Riesz lemma there is a unique bounded operator 
$A\in\Cal B(\Cal E)$ such that 
$$
\<Au,v\>_{\Cal E}=\rho(uv^*), \qquad u,v\in\Cal E.  
$$
$A$ is obviously a positive operator and in fact we 
have ${\text{trace}}\,A=1$, since for any 
orthonormal basis $v_1,v_2,\dots$ for $\Cal E$ 
$$
{\text{trace}}\,A=\sum_n\<Av_n,v_n\>=\sum_n\rho(v_nv_n^*)=
\rho(\alpha(\bold 1))=\rho(\bold 1)=1.
$$  
The following result shows how to compute the eigenvalue 
list of the restriction of $\rho$ to 
the commutant of $\alpha(\Cal B(H))$ 
in terms of this ``correlation" operator $A$.  

\proclaim{Lemma 4.3}
Let $\rho$ be a normal state of $\Cal B(H)$ and let 
$A$ be the positive trace class operator on $\Cal E$ defined by 
$\<Au,v\>_{\Cal E}=\rho(uv^*)$, $u,v\in\Cal E$.  Then 
$$
\Lambda(\rho\restriction_{\alpha(\Cal B(H))^\prime})=\Lambda(A).
$$ 
\endproclaim
\demo{proof}By Proposition 2.3.1,  
it suffices to exhibit a normal $*$-isomorphism 
$\theta$ of $\Cal B(\Cal E)$ onto $\alpha(\Cal B(H))^\prime$ 
with the property that
$$
\rho(\theta(T))={\text{trace}}(AT), \qquad T\in\Cal B(\Cal E).  
\tag{4.4}
$$
Consider the tensor product of Hilbert spaces $\Cal E\otimes H$.  
In order to define $\theta$ we claim first that there is 
a unique unitary operator $W: \Cal E\otimes H\to H$ which 
satisfies $W(v\otimes\xi)=v\xi$, $v\in\Cal E$, $\xi\in H$.  
Indeed, for $v,w\in\Cal E$, $\xi,\eta\in H$ we have 
$$
\<v\xi, w\eta\>_H = \<w^*v\xi,\eta\>=\<v,w\>_{\Cal E}\<\xi,\eta\>=
\<v\otimes\xi,w\otimes\eta\>_{\Cal E\otimes H}.
$$
It follows that 
there is a unique isometry $W: \Cal E\otimes H\to H$ with the 
stated property.  $W$ is unitary because its range spans 
all of $H$ (indeed, any vector $\zeta$ orthogonal to the range of 
$W$ has the property $v^*\zeta=0$ for every $v\in\Cal E$, hence
$\zeta=\alpha(\bold 1)\zeta=\sum_n v_nv_n^*\zeta=0$).

For every $X\in \Cal B(H)$ we have 
$$
W(\bold 1\otimes X)v\otimes \xi = W(v\otimes X\xi)=
vX\xi =\alpha(X)v\xi=\alpha(X)W(v\otimes\xi),  
$$
hence $W(\bold 1\otimes X)W^*=\alpha(X)$.  It follows
that 
$\alpha(\Cal B(H))^\prime=W(\Cal B(\Cal E)\otimes\bold 1)W^*$, 
and thus we can define a $*$-isomorphism 
$\theta : \Cal B(\Cal E)\to \alpha(\Cal B(H))^\prime$ by 
$\theta(T)=W(T\otimes\bold 1)W^*$.  

Writing $u\times\bar v$ for the rank-one operator on 
$\Cal E$ defined by $u\times\bar v:w\mapsto \<w,v\>_{\Cal E}u$, 
we claim that 
$$
\theta(u\times\bar v)=uv^*, \qquad {\text{for every }}
u,v\in\Cal E.
\tag{4.5}
$$
Indeed, if we pick a vector in $H$ of the form 
$\eta=w\xi=W(w\otimes\xi)$ where 
$w\in\Cal E$ and $\xi\in H$ then we have
$$
\align
\theta(u\times\bar v)\eta&=\theta(w\times\bar v)W(w\otimes\xi)
= W((u\times\bar v)\otimes\bold 1)w\otimes\xi=
W((u\times\bar v)w\otimes\xi)\\
&=\<w,v\>_{\Cal E}W(u\otimes\xi)=\<w,v\>_{\Cal E}u\xi=
uv^*w\xi=uv^*\eta,
\endalign
$$
and (4.5) follows because $H$ is spanned by all such vectors
$\eta$.   

Now for every 
rank-one operator $T=u\times\bar v\in\Cal B(E)$ we 
have 
$$
\rho(\theta(T))=\rho(\theta(u\times\bar v))=
\rho(uv^*)=\<Au,v\>_{\Cal E}=
{\text{trace}}(AT).  
$$
Formula (4.4) follows for finite rank $T\in\Cal B(\Cal E)$
by taking linear combinations, and the general case follows
by approximating an arbitrary operator $T\in\Cal B(\Cal E)$
in the strong operator topology 
with finite dimensional compressions $PTP$, $P$ ranging over 
an increasing sequence of finite dimensional projections 
with limit $\bold 1$.
\qedd\enddemo

The following formulas provide a key step.  

\proclaim{Lemma 4.6}
Let $\alpha$, $\Cal E$ be as above, let 
$\rho$ be a normal state of $\Cal B(H)$ and let 
$R\in\Cal L^1(H)$ be its density operator 
$\rho(X)={\text{trace}}(R X)$, $X\in\Cal B(H)$.  
Define a linear operator $L$ from $\Cal E$ into the 
Hilbert space $\Cal L^2(H)$ of all Hilbert-Schmidt
operators on $H$ by $Lv = R^{1/2}v$, $v\in\Cal E$.  
Then
\roster
\item"{4.6.1}"
$\<L^*Lu,v\>_{\Cal E}=\rho(uv^*), \qquad u,v\in\Cal E$, and 
\item"{4.6.2}"
For all $\xi_1,\xi_2,\eta_1,\eta_2\in H$ we have 
$$
\<LL^*(\xi_1\times\bar\xi_2),\eta_1\times\bar\eta_2\>_{\Cal L^2(H)}
= \<\alpha(\eta_2\times\bar\xi_2)R^{1/2}\xi_1,R^{1/2}\eta_1\>_H. 
$$
\endroster
\endproclaim

\demo{proof of (4.6.1)}
Simply write
$$
\<L^*Lu,v\>_{\Cal E}=\<Lu,Lv\>_{\Cal L^2(H)}=
\<R^{1/2}u,R^{1/2}v\>_{\Cal L^2(H)}=
{\text{trace}}(v^*Ru)=\rho(uv^*).  
$$

\demo{proof of (4.6.2)}
We have 
$$
\<LL^*(\xi_1\times\bar\xi_2),\eta_1\times\bar\eta_2\>_{\Cal L^2(H)}=
\<L^*(\xi_1\times\bar\xi_2),L^*(\eta_1\times\bar\eta_2)\>_{\Cal E}.  
\tag{4.7}
$$
Pick an orthonormal basis $v_1,v_2,\dots$ for $\Cal E$.  Then the 
right side of (4.7) can be rewritten as follows
$$
\align
&\sum_n\<L^*(\xi_1\times\bar\xi_2),v_n\>_{\Cal E}
\<v_n,L^*(\eta_1\times\bar\eta_2)\>_{\Cal E}=\\
&\sum_n\<\xi_1\times\bar\xi_2,R^{1/2}v_n\>_{\Cal L^2(H)}
\<R^{1/2}v_n,\eta_1\times\bar\eta_2\>_{\Cal L^2(H)}=\\
&\sum_n{\text{trace}}(v_n^*R^{1/2}\xi_1\times\bar\xi_2)
{\text{trace}}(R^{1/2}v_n\eta_2\times\bar\eta_1)=\\
&\sum_n\<v_n^*R^{1/2}\xi_1,\xi_2\>_H\<R^{1/2}v_n\eta_2,\eta_1\>_H.  
\endalign
$$

On the other hand, 
$$
\align
&\<\alpha(\eta_2\times\bar\xi_2)R^{1/2}\xi_1,R^{1/2}\eta_1\>_H=
\sum_n\<v_n(\eta_2\times\bar\xi_2)v_n^*R^{1/2}\xi_1,R^{1/2}\eta_1\>_H=\\
&\sum_n\<(\eta_2\times\bar\xi_2)v_n^*R^{1/2}\xi_1,v_n^*R^{1/2}\eta_1\>_H=
\sum_n\<v_n^*R^{1/2}\xi_1,\xi_2\>_H\<\eta_2,v_n^*R^{1/2}\eta_1\>_H,
\endalign
$$
and the last expression agrees with the bottom line of the 
previous formula.  
\qedd\enddemo

\proclaim{Lemma 4.8}
For a pair $A,B$ of self-adjoint compact operators on $H$, let
$A\circ B$ be the bounded operator defined on the Hilbert 
space $\Cal L^2(H)$ of Hilbert-Schmidt operators by 
$A\circ B(T)=ATB$.  Then $A\circ B$ is unitarily equivalent 
to $A\otimes B\in\Cal B(H\otimes H)$.  
\endproclaim
\demo{proof}
Pick orthonormal bases $e_1,e_2,\dots$ and $f_1,f_2,\dots$ for 
$H$ consisting of eigenvectors of $A$ and $B$, 
$Ae_n=\alpha_ne_n$, $Bf_n=\beta_nf_n$, $n=1,2,\dots$.  Letting
$e_m\times\bar f_n$ be the rank-one operator 
$\zeta\mapsto \<\zeta,f_n\>e_n$, then 
$\{e_m\times\bar f_n: m,n=1,2,\dots\}$ is an orthonormal basis 
for $\Cal L^2(H)$ and we have 
$$
A\circ B(e_m\times\bar f_n)=\alpha_m\beta_ne_m\times\bar f_n,
\qquad m,n=1,2,\dots.  
$$
Thus the unitary operator $W: \Cal L^2(H)\to H\otimes H$ 
defined by $W(e_m\times\bar f_n)=e_m\otimes f_n$, $m,n=1,2,\dots$
satisfies
$W (A\circ B)(e_m\times\bar f_n)=(A\otimes B )W(e_m\times\bar f_n)$
for every $m,n=1,2,\dots$ and hence $W( A\circ B) W^*=A\otimes B$.  
\qedd\enddemo

\demo{proof of Theorem C}
Let $R\in\Cal B(H)$ be the density operator of the normal 
state $\rho$, ${\text{trace}}(RT)=\rho(T)$, $T\in\Cal B(H)$. 
For every $t>0$ let $\Cal E_t$ be the Hilbert space of 
intertwining operators associated with $\alpha_t$,
$$
\Cal E_t=\{T\in\Cal B(H): \alpha_t(A)T=TA, \quad A\in\Cal B(H)\}, 
$$
and let $L_t: \Cal E_t\to \Cal L^2(H)$ be the operator of 
Lemma 3, $L_tv=R^{1/2}v$, $v\in\Cal E_t$.  

(4.6.1) implies that 
$\rho(uv^*)=\<L_t^*L_tu,v\>_{\Cal E}$, hence 
the correlation 
operator of $\rho\restriction_{\alpha_t(\Cal B(H))^\prime}$ 
is $L_t^*L_t$.  By Lemma 4.3 
$$
\Lambda(L_t^*L_t)=
\Lambda(\rho\restriction_{\alpha_t(\Cal B(H))^\prime}).  
$$
On the other hand, (2.1.3) implies that 
$\Lambda(L_t^*L_t)=\Lambda(L_tL_t^*)$.  Thus it suffices 
to show that the eigenvalue lists of the operators 
$L_tL_t^*\in\Cal B(\Cal L^2(H))$ converge to 
$\Lambda(\rho\otimes\omega)$, as $t\to\infty$, in the metric 
of eigenvalue lists.  

By (4.6.2) we have 
$$
\<L_tL_t^*(\xi_1\times\bar \xi_2),
\eta_1\times\bar\eta_2\>_{\Cal L^2(H)}=
\<\alpha_t(\eta_2\times\bar\xi_2)R^{1/2}\xi_1,R^{1/2}\eta_1\>_H, 
\tag{4.9}
$$
for all $\xi_1,\xi_2,\eta_1,\eta_2\in H$.  Now since 
$\alpha$ is pure, $\alpha_t(X)$ converges in the weak$^*$-topology
to $\omega(X)\bold 1$ as $t\to\infty$
(indeed, for every normal state $\sigma$, 
$\sigma(\alpha_t(X))$ converges to $\omega(X)=\sigma(\omega(X)\bold 1)$, 
and the assertion follows because every element of the predual of 
$\Cal B(H)$ is a linear combination of normal states).  
Thus if we take the limit on $t$ in the right side of (4.9) 
we obtain
$$
\align
\lim_{t\to\infty}
\<\alpha_t(\eta_2\times\bar\xi_2)R^{1/2}\xi_1,R^{1/2}\eta_1\>_H
&=\omega(\eta_2\times\bar\xi_2)\<R^{1/2}\xi_1,R^{1/2}\eta_1\>_H\\
=
\<\Omega\eta_2,\xi_2\>_H\<R\xi_1,\eta_1\>_H, 
\endalign
$$
where $\Omega$ is the density operator of $\omega$, 
$\omega(T)={\text{trace}}(\Omega T)$, $T\in\Cal B(H)$.  

Let $R\circ\Omega$ be the operator on $\Cal L^2(H)$ defined 
in Lemma 4.8, and notice that the right side of the preceding 
expression is 
$\<R\circ\Omega(\xi_1\times\bar\xi_2),
\eta_1\times\bar\eta_2\>_{\Cal L^2(H)}$.  Indeed, by definition 
of $R\circ\Omega$ we have 
$R\circ\Omega(\xi_1\times\bar\xi_2)=
R\xi_1\times\overline{\Omega\xi_2}$, 
and 
$$
\align
&\<R\xi_1\times\overline{\Omega\xi_2},\eta_1\times\bar\eta_2\>_{\Cal L^2(H)}
={\text{trace}}(\eta_2\times\bar\eta_1\cdot R\xi_1\times\overline{\Omega\xi_2})
=\\
&\<R\xi_1,\eta_1\>_H{\text{trace}}(\eta_2\times\overline{\omega\xi_2})
=\<R\xi_1,\eta_1\>_H\<\eta_2,\Omega\xi_2\>_H,   
\endalign
$$
which, as asserted, agrees with the right side of the previous 
expression.  

Thus we have shown that 
$$
\lim_{t\to\infty}\<L_tL_t^*(A),B\>_{\Cal L^2(H)}=
\<R\circ\Omega(A),B\>_{\Cal L^2(H)}
$$
for rank-one operators $A,B\in\Cal L^2(H)$.  Now Lemma 4.8 
implies that $R\circ\Omega$ is unitarily equivalent to 
$R\otimes\Omega\in\Cal B(H\otimes H)$, and hence $R\circ\Omega$
is a positive trace class operator for which 
$$
\Lambda(R\circ\Omega)=\Lambda(R\otimes\Omega)=\Lambda(\rho\otimes\omega).  
$$
On the other hand, Lemma 4.2 implies that 
$$
\lim_{t\to\infty}{\text{trace}}|L_tL_t^*-R\circ\Omega|=0.  
$$
By the inequality (2.3.2)
we conclude that 
$$
\limsup_{t\to\infty}\|\Lambda(L_tL_t^*)-\Lambda(R\circ\Omega)\|
\leq \lim_{t\to\infty}{\text{trace}}|L_tL_t^*-R\circ\Omega|=0.  
$$
We have already seen that 
$\Lambda(R\circ\Omega)=\Lambda(\rho\otimes\omega)$, and that 
$\Lambda(L_tL_t^*)=
\Lambda(\rho\restriction_{\alpha_t(\Cal B(H))^\prime})$.   
Thus Theorem C is proved.  \qedd
\enddemo

We now readily deduce the interaction
inequality.  

\proclaim{Theorem B}
Let $(U,M)$ be an interaction with past and 
future states $\omega_-$ and $\omega_+$, and let 
$\bar\omega_-$ and $\bar\omega_+$ be their natural 
extensions to the local \cstar\ $\Cal A$.  Then
$$
\|\bar\omega_- - \bar\omega_+\| \geq
\|\Lambda(\omega_-\otimes\omega_-)-
\Lambda(\omega_+\otimes\omega_+)\|.    
$$
\endproclaim
\demo{proof}
Fix $\epsilon>0$.  By Theorem C we can 
find $T>0$ large enough so 
that for all $t>T$ we have 
$$
\|\Lambda(\omega_+\restriction_{\Cal A_{[0,t]}})-
\Lambda(\omega_+\otimes\omega_+)\|\leq \epsilon
$$
as well as 
$$
\|\Lambda(\omega_-\restriction_{\Cal A_{[-t,0]}})-
\Lambda(\omega_-\otimes\omega_-)\|\leq \epsilon.  
$$

Now for $t\geq T$, 
$$
\align
\|\bar\omega_+ - \bar\omega_-\|&=
\|\bar\omega_+\circ\gamma_t - \bar\omega_-\circ\gamma_{-t}\|
\geq \|\bar\omega_+\circ\gamma_t\restriction_{\Cal A_{[-t,t]}}
-\bar\omega_-\circ\gamma_{-t}\restriction_{\Cal A_{[-t,t]}}\|
\\&=\|\omega_+\circ\gamma_t\restriction_{\Cal A_{[-t,t]}}
-\omega_-\circ\gamma_{-t}\restriction_{\Cal A_{[-t,t]}}\|.
\tag{4.10}
\endalign
$$
Since $\gamma_t$ gives rise to a 
$*$-isomorphism of $\Cal A_{[-t,t]}$
onto $\Cal A_{[0,2t]}$ while $\gamma_{-t}$ 
gives rise to a $*$-isomorphism of 
$\Cal A_{[-t,t]}$ onto $\Cal A_{[-2t,0]}$,   
(2.3.1) implies that 
$$
\align
\Lambda(\omega_+\circ\gamma_t\restriction_{\Cal A_{[-t,t]}})
&=\Lambda(\omega_+\restriction_{\Cal A_{[0,2t]}}), 
\qquad{\text{and}}\\
\Lambda(\omega_-\circ\gamma_{-t}\restriction_{\Cal A_{[-t,t]}})
&=\Lambda(\omega_-\restriction_{\Cal A_{[-2t,0]}}).
\endalign
$$  
Thus by Proposition 2.3 the last term of (4.10) 
is at at least
$$
\|\Lambda(\omega_+\restriction_{\Cal A_{[0,2t]}})-
\Lambda(\omega_-\restriction_{\Cal A_{[-2t,0]}})\|
$$
which by our initial choice of $T$ is 
at least 
$$
\|\Lambda(\omega_+\otimes\omega_+)-
\Lambda(\omega_-\otimes\omega_-)\| - 2\epsilon.  
$$
Since $\epsilon$ is arbitrary, the asserted inequality
follows.  
\qedd\enddemo

\proclaim{Corollary 1}
Let $(U,M)$ be an interaction with past and future states
$\omega_-$, $\omega_+$.  If $\Lambda(\omega_-)\neq\Lambda(\omega_+)$, 
then the interaction is nontrivial.  
\endproclaim
\demo{proof}
Contrapositively, 
suppose that the interaction is trivial and let  
$\Omega_-$ and $\Omega_+$ be the 
respective density operators of $\omega_-$ and 
$\omega_+$.  Theorem B 
implies that $\Omega_-\otimes\Omega_-$ and 
$\Omega_+\otimes\Omega_+$ must have the same eigenvalue
list.   (2.1.4) implies that for every 
$n=1,2,\dots$ we have
$$
{\text{trace}}(\Omega_-^n)^2=
{\text{trace}}((\Omega_-\otimes\Omega_-)^n)=
{\text{trace}}((\Omega_+\otimes\Omega_+)^n)=
{\text{trace}}(\Omega_+^n)^2.  
$$ 
Taking the square root we find that  
${\text{trace}}(\Omega_-^n)={\text{trace}}(\Omega_+^n)$ 
for every $n=1,2,\dots$ and another application of 
(2.1.4) leads to 
$\Lambda(\Omega_-)=\Lambda(\Omega_+)$.  
\qedd\enddemo

\proclaim{Corollary 2}Let $n=1,2,\dots,\infty$ and 
choose $\epsilon>0$.  There is an interaction 
$(U,M)$ whose past and future $E_0$-semigroups 
are cocycle-conjugate to the $CAR/CCR$ flow of 
index $n$ such that 
$$
\|\bar\omega_+ - \bar\omega_-\|\geq 2-\epsilon.  
$$
\endproclaim
\demo{proof}
Choose positive integers $p<q$ and consider the eigenvalue 
lists 
$$
\align
\Lambda_-&=\{1/p,1/p,\dots,1/p,0,0,\dots\}\\
\Lambda_+&=\{1/q,1/q,\dots,1/q,0,0,\dots\},  
\endalign
$$
where $1/p$ is repeated $p$ times and $1/q$ is 
repeated $q$ times.  
Theorem A implies that there is an interaction $(U,M)$ whose
past and future $E_0$ semigroups are cocycle-conjugate 
to the $CAR/CCR$ flow of index $n$, for which 
$\Lambda(\omega_-)=\Lambda_-$
and $\Lambda(\omega_+)=\Lambda_+$.  By 
Theorem B 
$$
\|\bar\omega_+-\bar\omega_-\|\geq 
\|\Lambda(\omega_+\otimes\omega_+)-
\Lambda(\omega_-\otimes\omega_-)\|.  
$$
If we neglect zeros, 
the eigenvalue list of $\omega_-\otimes\omega_-$ 
consists of the single eigvalue $1/p^2$, repeated 
$p^2$ times, and that of $\omega_+\otimes\omega_+$ consists
of $1/q^2$ repeated $q^2$ times.  Thus 
$$
\|\Lambda(\omega_+\otimes\omega_+)-
\Lambda(\omega_-\otimes\omega_-)\|=
p^2(1/p^2-1/q^2)+(q^2-p^2)/q^2=2-2p^2/q^2,
$$
and the inequality of Corollary 2 follows whenever  
$q$ is larger than $p\sqrt{2/\epsilon}$.  
\qedd\enddemo

\end